\documentclass{amsart}
\usepackage{amsmath,amssymb,amscd}

\newcommand\al{\alpha}
\newcommand\be{\beta}
\newcommand\eps{\varepsilon}
\newcommand\la{\lambda}

\newcommand\BA{\mathbb A}

\newcommand\CX{{\mathcal X}}
\newcommand\CY{{\mathcal Y}}

\newcommand\gm{\mathfrak m}

\newcommand\gp{{\mathfrak p}}
\newcommand\gq{{\mathfrak q}}

 \DeclareMathOperator\Spec{Spec}

 \DeclareMathOperator\Fr{Fr}

\newcommand\isom{\simeq}

\newcommand\chr{\operatorname{char}}
\newcommand\Gal{\operatorname{Gal}}

 \DeclareMathOperator\su{{su}}

\theoremstyle{plain} \swapnumbers

\newtheorem{thm}[subsection]{Theorem}
\newtheorem{lm}[subsection]{Lemma}
\newtheorem{cor}{Corollary}[subsection]

\theoremstyle{definition}

\newtheorem{rk}[cor]{Remark}

\newcommand\paper[1]{\emph{#1}, }
\newcommand\jour[1]{#1 }
\newcommand\vol[1]{{\bf#1}}
\newcommand\yr[1]{ (#1)}

\begin{document}

\copyrightinfo{2002}
    {American Mathematical Society}

\author{Igor Zhukov}

\title{Ramification of Surfaces: \\ Artin-Schreier Extensions}
\address{Department of Mathematics and Mechanics, St. Petersburg University \\
Bibliotechnaya pl., 2, Staryj Petergof \\
198904 St. Petersburg, Russia
}

\email{Igor.Zhukov@mail.ru}

\thanks{I acknowledge the financial support from
Royal Society and from RFBR (projects 00-01-00140 and
01-01-00997.)}

\dedicatory{Dedicated to A.~N.~Parshin on occasion of his 60th birthday}



\date{April 30, 2002}

\begin{abstract}
Let $A$ be a regular 2-dimensional local ring
of characteristic $p>0$, and let $L/K$
be a cyclic extension of degree $p$ of its field
of fractions such that the corresponding branch divisor
is normal crossing. For each $\gp\in\Spec A$ of height 1
such that $A/\gp$ is regular,
consider the ramification jump $h_\gp$ of the
extension of the residue field at $\gp$.
In this paper the semi-continuity of $h_\gp$
with respect to Zariski topology of suitable jet
spaces is proved. The asymptotic of $h_\gp$ with
respect to intersection multiplicity of the prime divisor
defined by $\gp$
and the branch divisor is also addressed.
\end{abstract}

\maketitle


Let $L/K$ be a finite Galois extension of
the function field of a connected normal $n$-dimensional scheme $\CX$,
and let $f:\CY\to\CX$ be the normalization of $\CX$ in $L$.
Denote by $R\subset\CX$ the branch locus of this morphism.
One can attach numerous ramification invariants to the components
of $R$: different and discriminant, depth of ramification, lower and upper ramification
subgroups, genome, generalized Artin and Swan conductors etc.
If $n=1$, all these invariants are related with each other in a
very nice way. In particular, all the invariants can be computed
in terms of (lower) ramification filtrations. However, if $n\ge2$,
the relations between these invariants are, as a rule, obscure. It would
be desirable to construct some basic system of ramification
invariants of $(\CX,L/K)$ that determines all the other invariants.

In the present paper and in \cite{ram:su}
we start to develop an approach to ramification theory
of surfaces which can be described roughly as follows.

Let $C$ be a curve on a 2-dimensional scheme $\CX$, and let $D$ be any curve
on $\CY$ over $C$. Then the
natural morphism $D\to C$ has well known ramification invariants,
namely, the (lower) ramification filtrations at all points where
$C$ meets $R$.
The idea is to collect these invariants for all regular
curves on $C$ that are distinct from the components of $R$ and to consider
these data as a system of invariants of $(\CX,L/K)$.
(Some properties of these invariants were discussed in \cite{Del:l} and
\cite{Br:sur}.)
In \cite{ram:su} a more detailed description of this program is given.

In this paper we treat only the case when $\CX=\Spec A$,
$A$ is a regular 2-dimensional local ring of characteristic $p>0$,
and $L/K$ is a cyclic extension of degree $p$. Then the
above-mentioned ramification invariants are reduced to the
set of numbers $h_\gp(L/K)$, where $\gp$ runs over the set of
prime ideals of $A$ such that $A/\gp$ is a 1-dimensional regular
ring, and $L/K$ is unramified at the place $\gp$. Here
$h_\gp(L/K)$ denotes the only ramification jump of the extension
of residue fields emerging at $\gp$, if this extension is
non-trivial; otherwise, $h_\gp(L/K)=0$.

Theorems \ref{usu}--\ref{asymp} describe the behavior of
$h_\gp(L/K)$ as one varies $\gp$. These theorems give positive
answers to a part of the questions in \cite[\S2]{ram:su} for $\CX$
and $L/K$ as above.

I am very much grateful to A. N. Parshin for his permanent
encouragement in my attempts to analyze ramification invariants
of higher-dimensional schemes. I am also very much indebted
to V. P. Snaith for the invitation to Southampton and for a lot of
ramification discussions there and in St. Petersburg.

\section{Terminology and notation}

Let $A$ be an equal characteristic regular 2-dimensional local ring
with the maximal ideal $\gm=(T,U)$ and perfect residue field $k$,
$N$ a positive integer.
We have a canonical isomorphism between the completion of $A$
and $k[[T,U]]$.
A map $k=A/\gm\to A$ is said to
be a section of level $N$ if the diagram
$$
\begin{CD}
k & @>\lambda>> & A
\\
@VVV && @VVV
\\
k[[T,U]]/(T,U)^{N}  & @= & A/\gm^{N}
\end{CD}
$$
commutes.

For any commutative ring $A$ denote
$$
\Spec_1A=\{\gp\in\Spec A|\dim(A/\gp)=1\}.
$$

For a prime ideal $\gp\in\Spec_1A$,
denote by $F_\gp$ the prime divisor $\Spec(A/\gp)$ in
$\Spec A$.

If $R$ is a domain, $\Fr R$ is the field of fractions of $R$.

If $A$ is a local ring, $\widehat{A}$ is the completion of $A$.

\section{Theorems}
\label{Theorems}

Let $A$ be a regular two-dimensional local ring, $K=\Fr A$,
$\gm$ the maximal ideal of $A$, $k$ the
residue field assumed to be algebraically closed, $\chr k=p>0$.

For any two distinct prime divisors $F_\gp$, $F_{\gp'}$
we define their intersection number as
$$
(F_\gp.F_{\gp'})=\dim_k A/(\gp+\gp');
$$
by linearity this definition can be extended to any two
divisors $C,C'$ with no common components.

Let $L/K$ be a cyclic extension of degree $p$, and let $B$ be the
integral closure of $A$ in $L$. Let $\gp_1,\dots,\gp_d\in\Spec_1A$
be those prime ideals where $L/K$ is ramified. In this paper we
shall assume that the branch divisor is (strictly) normal crossing, i.~e.,
$d\le2$; $A/\gp_1,\dots,A/\gp_d$ are regular; if $d=2$, then
$(F_{\gp_1}.F_{\gp_2})=1$. We are not interested in unramified
extensions; thus, we shall assume $d\ge1$. The valuations on $K$
associated with prime divisors $F_{\gp_1},\dots F_{\gp_d}$ will be
denoted by $v_1,\dots,v_d$. Next,
\[
 U_A=\{\gp\in\Spec_1 A|A/\gp\text{ is regular,
 }\gp\ne\gp_1,\dots,\gp_d\}.
\]

Let $\gp\in U_A$.
Consider $D_\gq$, the decomposition subgroup in $\Gal(L/K)$
at $\gq$, where $\gq$ is a prime ideal of $B$ over $\gp$.
Since $L/K$ was not ramified at
$\gp$, we have $D_\gq\isom\Gal(L(\gq)/K(\gp))$, where $K(\gp)$ is the
fraction field of $A/\gp$, and $L(\gq)$ is the fraction field of $B/\gq$
(see \cite[Ch. V, \S2, Prop. 6]{Bourbaki}).

The ring $A/\gp$ is a discrete valuation
ring with residue field $k$,
and we obtain (lower) ramification filtration
on the group $\Gal(L(\gq)/K(\gp))$. The  ramification
subgroups will be denoted by  $G_{i,\gp}$, $i\ge-1$. We put
\[
h_{\gp}(L/K)= \begin{cases}
\max\{i|p=|G_{i,\gp}|\}, & e(L(\gq)/K(\gp))=p,
\\
0, & \text{otherwise}.
\end{cases}
\]
It is easy to see that $h_{\gp}(L/K)$ is well defined.
(To relate this notation with that from \cite{ram:su},
denote by $O$ the closed point of $\Spec A$. Then
$h_{\gp}(L/K)=w_{F_\gp,O}^{(1)}(L/K)$.)

For $\gp\in U_A$,
the set
$$
J_n(\gp)=\{\gp'\in U_A|(F_{\gp'}.F_\gp)\ge n+1\}
$$
is said to be the jet of $\gp$ of order $n$. We introduce also
$$
T_r= \begin{cases}
\{\gp\in U_A|(F_\gp.F_{\gp_1})=r\}, & d=1,
\\
\{\gp\in U_A|(F_\gp.F_{\gp_1})=r,\,(F_\gp.F_{\gp_2})=1\}, & d=2,
\end{cases}
$$
and
$$
T_{r,n}=\{J_n(\gp)|\gp\in T_r\}.
$$

\begin{thm}
\label{usu}
\emph{(existence of a uniform sufficient jet order)}
For any $r\ge1$ there
exists $s$ such that if $\gp,\gp'\in T_r$ and $(F_\gp.F_{\gp'})\ge
s+1$, then $h_\gp(L/K)=h_{\gp'}(L/K)$. Let $\su_{1,r}(L/K)$
be the minimal such $s$. Then there exists
$N\ge1$ such that $\su_{1,r}(L/K)< N r$ for any $r$.
\end{thm}

This theorem is a particular case of Corollary 4.1.1 in
\cite{ram:su}. Moreover, if $L=K(x)$, $x^p-x=a\in K$,
and $m_i=-v_i(a)$, $i=1,\dots,d$, it follows from
\cite[Prop. 6.1]{ram:su} that
\[
\su_{1,r}(L/K)+1\le
\begin{cases} (m_1+1)r+m_2+1, &d=2, \\ (m_1+1)r, &d=1. \end{cases}
\]
Thus, one can take $N=\sum_{i=1}^d (m_i+1)$.

\begin{rk}
In fact, we need not refer to \cite{ram:su} since both Theorem
\ref{usu} and the upper bound for $\su_{1,r}(L/K)$ follow from the
argument in subsection \ref{pf1}.
\end{rk}

To state further theorems, we have to introduce a structure
of affine space on $T_r$.
If $d=2$, denote by $T$ and $U$ any generators
of $\gp_1$ and $\gp_2$ respectively. If $d=1$, denote by
$(T,U)$ any system of regular local parameters in $A$ such
that $T$ is a generator of $\gp_1$. Then $\widehat{A}$ can be
identified with $k[[T,U]]$.

Fix a positive integer $n$ and a section $\lambda:A/\gm\to A$ of level $n$.

First, we consider the case $r=1$, $d=1$.

Let $\gp\in T_1$ (i.~e., this is the ideal of the germ of a curve
transversal to $F_{\gp_1}$). Then $\gp=(f)$, where $f\not\equiv0\mod(T,U^2)$,
and we may assume $f\equiv -U\mod(T,U^2)$.
By Weierstra\ss\ preparation theorem
\cite[Ch. VII, \S3, Prop. 6]{Bourbaki} there exists a unique
$\eps\in(\widehat{A})^*$ such that
$f\eps=-U+\sum_{i=1}^\infty\al_iT^i$
for some $\al_i\in k$, $i=1,2,\dots$
Take $\eps_0\in A^*$ such that $\eps_0\equiv\eps\mod(T,U)^{n+1}$.
Replacing $f$ with $f\eps_0$, we may assume
without loss of generality that
\begin{equation}
f\equiv-U+\la(\al_1)T+\dots+\la(\al_n)T^n\mod\deg n+1,
\label{WA}
\end{equation}
where $\al_1,\dots,\al_n\in k$ are determined uniquely by $\gp$.
Thus, we can identify $T_{1,n}$ with the set of closed points of $\BA_k^n$
via
\begin{equation*}
(\al_1,\dots,\al_n) \mapsto J_n((-U+\la(\al_1)T+\dots+\la(\al_n)T^n)).
\end{equation*}
Note that $\al_1,\al_2,\dots$ are independent of $\lambda$ and
$n$. They are in fact the coefficients in the expansion
\[
u=\al_1t+\al_2t^2+\dots,
\]
where $t$ and $u$ are the images of
$T$ and $U$ in $\widehat{A/\gp}\isom k[[t]]$.

Next, let either $r\ge2$ or $d=2$, $r=1$.
If $\gp\in T_r$, we have $\gp=(f)$,
\begin{equation}
f\equiv-T+\la(\be_r)U^r+\dots+\la(\be_n)U^n\mod\deg{n+1},
\label{WB}
\end{equation}
where $\be_r,\dots,\be_n\in k$ are uniquely determined by $\gp$,
and $\be_r\ne0$.
We have a bijection
\begin{align*}
(\BA_k^{n+1-r})_{x_0\ne0} & \to T_{r,n},
\\
(\be_r,\dots,\be_n) & \mapsto J_n((-T+\la(\be_r)U^r+\dots+\la(\be_n)U^n)).
\end{align*}
Here $(\BA_k^{n+1-r})_{x_0\ne0}$ is the set of closed points of
$\BA_k^{n+1-r}$ with a non-vanishing first coordinate.
As in the previous case, $\be_r,\be_{r+1},\dots$ are independent of $\lambda$ and $n$.


All the below theorems will be proved in the next section.

\begin{thm}
\label{semicont}
\emph{(semi-continuity of a jump)}
Let $n\ge\su_{1,r}(L/K)$. Then for any $s\ge0$ the set
\[
\{J_n(\gp)|\gp\in T_r;\,h_\gp(L/K)\le s\}
\]
is a closed subset in $T_{r,n}$.
\end{thm}

\begin{thm}
\label{gener}
\emph{(generic value of a jump)}
The supremum
$$
h_r(L/K)=\sup\{h_\gp(L/K)|\gp\in T_r\}
$$
is finite.
\end{thm}

Recall that a local ring $A$ is a G-ring
iff for any $\gp\in\Spec A$ and any finite extension
$L/\Fr(A/\gp)$ the formal fiber $\widehat{A}\otimes_AL$
is a regular ring, see \cite[(33.C)]{Matsumura}.
(This is one of the conditions in the definition
of an excellent ring. Local rings of varieties are G-rings.)

\begin{thm}
\label{asymp}
\emph{(asymptotic of jumps)}
Assume in addition that $A$ is a G-ring.
The sequence $(h_r(L/K)/r)_r$ is convergent.
\end{thm}

These three theorems supply positive answers to the questions with
the same names in \cite{ram:su} for the case when $\CX$ is local,
and $L/K$ is a cyclic extension of degree $p$. In the question
related to the asymptotic of jumps one needs an additional
assumption that $A$ is a G-ring.

\section{Proofs}

\begin{lm}
\label{approx}
Let $A$ be a regular 2-dimensional local ring,
$K$ the fraction field of $A$,
$\gp_1,\dots,\gp_n,\gq$ distinct primes of height 1 in $A$,
$x_1,\dots,x_n\in K$.
Then there exists $x\in K$ such that $x-x_i\in A_{\gp_i}$, $i=1,\dots,n$,
and $x\in A_\gp$, $\gp\ne\gp_1,\dots,\gp_n,\gq$.
\end{lm}


\begin{proof}
Reduction to the case $n=1$ is immediate.
Since $A$ is regular,
$\gp=\gp_1=(t)$ and $\gq=(u)$ for some $t,u\in A$.
One can write $x_1=t^{-r}a$, $a\in A_\gp$, $r$ is integral.
We may assume $r>0$; otherwise, take $x=0$.

Apply induction on $r$; let $r=1$. Observe that $u+\gp$
is a non-zero element of the maximal ideal of a (not necessarily
regular) one-dimensional local ring $A/\gp$. Therefore,
$A_\gp/\gp A_\gp=(A/\gp)_u$. Then $au^l\equiv a_0\mod\gp A_\gp$, where
$a_0\in A$, $l$ is a non-negative integer. One can take $x=t^{-1}u^{-l}a_0$.

Let $r>1$. Applying the case $r=1$, one can find $x'\in A[t^{-1},u^{-1}]$ such that
$x'-t^{r-1}x_1\in A_\gp$. It remains to apply the induction
hypothesis, taking $x_1-t^{-(r-1)}x'$ for a new $x_1$.
\end{proof}

The following lemma is standard and easy to prove.

\begin{lm}
\label{Hyodo}
Let $R$ be a discrete valuation ring of prime characteristic $p$,
$K$ the fraction field,
$\pi$ a uniformizer,
$L/K$ a  cyclic extension of degree $p$ which is not unramified. Then $L=K(x)$,
$x^p-x=\pi^{-m}a$, $m>0$, $a\in R^*$, and one of the
following two conditions holds:

(i) $(m,p)=1$, and $L/K$ is totally ramified;

(ii) $p|m$, $a\notin R^p+(\pi)$, and $L/K$ is ferociously
ramified, i.~e., the inseparable degree of the extension of
residue fields is equal to $[L:K]$.
\end{lm}

\begin{cor}
\label{AS}
Let $A$ be a regular 2-dimensional local ring of characteristic $p>0$, $L/K$
a cyclic extension of the fraction field of degree $p$.

1. Assume that the branch locus of $L/K$ in $\Spec A$ is regular,
i.~e., it consists of a single regular prime divisor $F_{\gq_0}$;
we have $\gq_0=(t)$ for some $t\in A$. Then $L=K(x)$, $x^p-x=a$,
$a\in A[t^{-1}]$.

2.  Assume that the branch locus of $L/K$ in $\Spec A$ consists of
two transversal regular prime divisors $F_{\gq_0}$, $F_{\gq_1}$,
where $\gq_0=(t)$ and $\gq_1=(u)$. Then $L=K(x)$, $x^p-x=a$, $a\in
A[t^{-1},u^{-1}]$.
\end{cor}

\begin{proof}
By Artin-Schreier theory, $L=K(x_0)$, $x_0^p-x_0=a_0$, $a_0\in K$.

By Lemma \ref{Hyodo}, in the case ``1'' (resp., ``2'') $\gq_i$ are poles of
$a_0$, where $i=0$ (resp., $i=0,1$).
Let $\gp_1,\dots,\gp_n\in\Spec A$ be all the other poles of $a_0$.
In the case ``1'' choose $u\in A$ such that $(t,u)$
are local parameters of $A$ and $\gq=(u)$ is distinct from
$\gp_1,\dots,\gp_n$. In the case ``2'' put $\gq=\gq_1$.

Since $L/K$ is not ramified
at places $\gp_1,\dots,\gp_n$, by Lemma \ref{Hyodo} there exist $d_i\in K$
such that $a_0+d_i^p-d_i\in A_{\gp_i}$, $i=1,\dots,n$. By Lemma \ref{approx},
there exists $d\in K$ such that $d-d_i\in A_{\gp_i}$, $i=1,\dots,n$,
and $d$ is integral outside $\{\gp_1,\dots,\gp_n,\gq\}$.

Put $a_1=a_0+d^p-d$. Then $L=K(x_1)$, $x_1^p-x_1=a_1$,
and the only poles of $a_1$ are $\gp_0$ and (possibly) $\gq$.
It follows $a_1\in A[t^{-1},u^{-1}]$, and in the case ``2'' we are
done.

Finally, denote by $a$ any element of $A[t^{-1},u^{-1}]$ which
takes the form $a=a_1+D^p-D$, $D\in K$, and is of minimal possible
degree $r$ in $u^{-1}$. It remains to prove that $r=0$.

Assume $r>0$.
Applying Lemma \ref{Hyodo} to the completion $R$
of $A_\gq$, we see that $r=ps$ with integral $s$.
(Otherwise, $L/K$ would be totally ramified at $\gq$!)
Hence $a=u^{-ps}b$, where $b\in A[t^{-1}]$.
We have $b\in R^p+uR=A_\gq^p+uR=A[t^{-1}]^p+uR$.
(Otherwise, by Lemma \ref{Hyodo},
$L/K$ would be ferociously ramified at $\gq$!)
If
\[
b\equiv b_0^p\mod uR,\quad b_0\in A[t^{-1}],
\]
then in fact
$$
b-b_0^p\in A[t^{-1}]\cap uR = A[t^{-1}]\cap uA_\gq = uA[t^{-1}].
$$
Therefore,
$a'=a-(u^{-s}b_0)^p+u^{-s}b_0$ is in $u^{-(ps-1)}A[t^{-1}]$,
a contradiction with the choice of $a$.
\end{proof}

The lemma below is also standard and easy to prove.

\begin{lm}
\label{j}
In the case $(m,p)=1$ in Lemma \ref{Hyodo} the only ramification
jump of $L/K$ is $m$.
\end{lm}

\begin{lm}
\label{ej}
Let $R$ be a complete discrete valuation ring of prime
characteristic $p$ with algebraically closed residue field,
$K$ the field of fractions, $L_i=K(x_i)$, $x_i^p-x_i=a_i\in K$,
$i=1,2$. Assume that $a_1-a_2\in R$. Then $L_1=L_2$.
\end{lm}

\begin{proof}
We have $L_1L_2=L_1(x)$, $x^p-x=a_1-a_2$. Since the residue
field is algebraically closed, the polynomial $X^p-X-(a_1-a_2)$
has a root in the ring of integers of $L_1$, whence $L_1L_2=L_1$.
Similarly, $L_1L_2=L_2$.
\end{proof}

\subsection{Proof of Theorems \ref{gener} and \ref{semicont}}
\label{pf1}
\

Let $A,K,\gm,k,p$; $L/K$; $\gp_1,\dots,\gp_d$; $v_1,\dots,v_d$; $T,U$ be as in
section \ref{Theorems}.

By Corollary \ref{AS}, $L=K(x)$,
where $x^p-x=a$, $a\in A[T^{-1},U^{-1}]$ if $d=2$,
and $a\in A[T^{-1}]$ if $d=1$.

Denote $m=-v_1(a)$. Put $n=-v_2(a)$ if $d=2$, and $n=0$ otherwise.
Fix some $r\ge1$ and a section $\lambda:k=A/\gm\to A$ of level $rm+n+1$.
One can write
\begin{equation}
a=
T^{-m}U^{-n}\Bigl(\sum_{i=0}^{m+n-1}\sum_{j=0}^{rm+n-1}
               \lambda(\theta_{ij})T^iU^j+Z\Bigr),
\label{dth}
\end{equation}
where $\theta_{ij}\in k$, $Z\in T^{m+n}A+U^{rm+n}A$.
Looking at the expansion of $T^mU^na$ in the completion of $A$,
we conclude that $\theta_{ij}$ are independent of $r$ and $\lambda$.

We start with the case, when $r=1$ and $d=1$.
Let $\gp\in T_1$. Denote by $t$ and $u$ the images of $T$ and $U$
in $A/\gp$; $t$ is a local parameter in this discrete valuation ring.
Introducing $\al_1,\dots,\al_{m}$ as in \eqref{WA}, we get in the completion
of $A/\gp$:
$$
u\equiv\al_1t+\dots+\al_{m}t^{m}\mod t^{m+1}.
$$

By the definition,
$h_\gp(L/K)$ is the
only ramification jump of $k((t))(y)/k((t))$, where
\begin{equation}
\label{C}
\begin{split}
y^p-y &=
t^{-m}\Bigl(\sum_{i=0}^{m-1}\sum_{j=0}^{m-1}
               \theta_{ij}t^iu^j\Bigr)
               \\
& \equiv t^{-m}\Bigl(\sum_{i=0}^{m-1}\sum_{j=0}^{m-1}
               \theta_{ij}t^{i+j}(\al_1+\dots+\al_{m}t^{m-1})^{j}\Bigr)
               \mod k[[t]],
\end{split}
\end{equation}
if this extension is non-trivial, and 0 otherwise.

Now Theorem \ref{gener} for $r=1$ and $d=1$ follows from Lemma \ref{j}
and Lemma \ref{ej}; we have $h_1(L/K)\le m$.

Next, \eqref{C} implies
$$
y^p-y \equiv t^{-m}\sum_{i=0}^{m-1}F_i(\al_1,\dots,\al_i)t^i\mod k[[t]],
$$
where $F_0=\theta_{00}$, and $F_i\in k[X_1,\dots,X_i]$ are
determined by
\begin{displaymath}
\sum_{i=0}^{m-1}\sum_{j=0}^{m-1}
               \theta_{ij}\tau^{i+j}(X_1+\dots+X_{m}\tau^{m-1})^{j}\equiv
\sum_{i=0}^{m-1}F_i(X_1,\dots,X_i)\tau^i \mod \tau^m
\end{displaymath}
in $k[X_1,\dots,X_m,\tau]$.
Lemma \ref{ej} implies $k((t))(y)=k((t))(\tilde y)$,
$$
\tilde y^p-\tilde y = \sum_{\substack{ 1\le l\le m\\(p,l)=1}}
G_l(\al_1,\dots,\al_{m-1})t^{-l},
$$
where
$$
G_l(\al_1,\dots,\al_{m-1})=F_{m-l}(\al_1,\dots,\al_{m-l})
+F_{m-pl}(\al_1,\dots,\al_{m-pl})^{p^{-1}}+\dots
$$
For $s=0,\dots,m-1$, we have
$$
h_\gp(L/K)\le s
$$
iff
$$
G_l(\al_1,\dots,\al_{m-1})=0,\quad s+1\le l\le m, \, (p,l)=1.
$$
Therefore, for any $M\ge m-1$ the set
$$
\{J_M(\gp)|\gp\in T_1;\,h_\gp(L/K)\le s\}
$$
can be identified with the set of common zeroes of polynomials
$$
\sum_{\nu=0}^{N_l}F_{m-p^\nu l}^{p^{N_l-\nu}},\quad s+1\le l\le m, \,
(p,l)=1,
$$
where $N_l=\max\{i|m-p^il\ge0\}$, and all $F_i$ ($0\le i\le m-1$)
are considered as elements of $k[X_1,\dots,X_M]$.
This proves Theorem \ref{semicont} for $r=1$ and $d=1$.

Next, assume that either $r\ge2$ or $r=1$ and $d=2$.

Let $\gp\in T_r$. Denote by $t$ and $u$ the images of $T$ and $U$
in $A/\gp$; $u$ is a local parameter in this discrete valuation ring.
Introducing $\be_r,\dots,\be_{r+rm+n-1}$ as in \eqref{WB} (even if $r=1$),
we get in the completion of $A/\gp$:
$$
t\equiv\be_ru^r+\dots+\be_{r+rm+n-1}u^{r+rm+n-1}\mod u^{r+rm+n},
$$
and $\be_r\ne0$. The image of $T^{-m}U^{-n}Z$ in
$\widehat{A_\gp/\gp A_\gp}=k((u))$ belongs to $k[[u]]$.
Therefore, in view of Lemma \ref{ej}, $h_\gp(L/K)$ is the
only ramification jump of the extension $k((u))(y)/k((u))$, where
\begin{multline}
\label{CB}
y^p-y =
t^{-m}u^{-n}\Bigl(\sum_{i=0}^{m+n-1}\sum_{j=0}^{rm+n-1}
               \theta_{ij}t^iu^j\Bigr)
\equiv \be_r^{-m}u^{-rm-n} \times
\\
\times \Bigl(\sum_{i=0}^{m+n-1}\sum_{j=0}^{rm+n-1}
               \theta_{ij}
               \be_r^i
               u^{ri+j}
               \Bigl(\sum_{\nu=0}^{rm+n-1}
                     \frac{\be_{r+\nu}}{\be_r}u^\nu
                     \Bigr)^{-m+i}
               \Bigr)
               \mod k[[u]].
\end{multline}
Theorem \ref{gener} follows immediately; we have
$h_r(L/K)\le rm+n$.

It follows from \eqref{CB} that
\[
y^p-y \equiv \be_r^{-m}u^{-rm-n}
        \sum_{i=0}^{rm+n-1}
        F_i^{(r)}(\be_r,\be_{r+1}/\be_r,\dots,\be_{r+i}/\be_r)u^i
        \mod k[[u]],
\]
where $F_0^{(r)}=\theta_{00}$, and $F_i^{(r)}\in k[X_0,\dots,X_i]$
are determined by
\begin{multline}
\label{Fr}
\sum_{i=0}^{m+n-1}\sum_{j=0}^{rm+n-1}
               \theta_{ij}
               X_0^i
               \tau^{ri+j}
               \Bigl(1+\sum_{\nu=1}^{rm+n-1}
                     X_\nu\tau^\nu
                     \Bigr)^{-m+i}
    \equiv \\
    \equiv \sum_{i=0}^{rm+n-1}
        F_i^{(r)}(X_0,X_1,\dots,X_i)\tau^i \mod \tau^{rm+n}
\end{multline}
in $k[X_0,\dots,X_{rm+n-1}][[\tau]]$.

Lemma \ref{ej} implies
$k((u))(y)=k((u))(\tilde y)$,
$$
\tilde y^p-\tilde y = \sum_{\substack{ 1\le l\le rm+n\\(p,l)=1}}
G_l^{(r)}(\be_r,\dots,\be_{r+rm+n-1})u^{-l},
$$
where
$$
G_l^{(r)}(\be_r,\dots,\be_{r+rm+n-1})=
\sum_{\nu=0}^{N_l}
\Bigl(\be_r^{-m}
 F_{i(\nu)}^{(r)}\Bigl(\be_r,\frac{\be_{r+1}}{\be_r},\dots,\frac{\be_{r+i(\nu)}}{\be_r}\Bigr)
 \Bigr)^{p^{-\nu}},
$$
$i(\nu)=rm+n-p^\nu l$, $N_l=\max\{\nu|i(\nu)\ge0\}$. For
$s=0,\dots,rm+n-1$, we have
$$
h_\gp(L/K)\le s
$$
iff
$$
G_l^{(r)}(\be_r,\dots,\be_{r+rm+n-1})=0,\quad
s+1\le l\le rm+n, \, (p,l)=1.
$$
Therefore, for any $M\ge r+rm+n-1$ the set
$$\{J_M(\gp)|\gp\in T_r;\,h_\gp(L/K)\le s\}$$
can be identified with the set of all common zeroes of polynomials
$$
\left\{
\sum_{\nu=0}^{N_l}
X_0^NF_{rm+n-p^\nu l}^{(r)}
\Bigl(X_0,\frac{X_1}{X_0},\dots,\frac{X_{rm+n-p^\nu l}}{X_0}\Bigr)^{p^{N_l-\nu}}
\right\}_{{\quad s+1\le l\le rm+n, \, (p,l)=1}}
$$
in $(\BA_k^{M+1-r})_{x_0\ne0}$, where $N$ is big enough.
This proves Theorem \ref{semicont}.
\qed

\subsection{Proof of Theorem \ref{asymp}}

 Choose an Artin-Schreier equation
$x^p-x=a$ for $L/K$ with $a\in A[T^{-1},U^{-1}]$
such that $m=-v_1(a)$ is minimal, and define $\theta_{ij}$ as
in \eqref{dth}. Put $n=-v_2(a)$ if $d=2$,
and $n=0$ if $d=1$.

Case 1: $p\nmid m$.
Take the minimal $j$ such that $\theta_{0j}\ne0$; such $j$ exists
in view of minimality of $m$.
Take any $r>\max(1,j)$. First, let $p\nmid rm+n-j$.
Since $\theta_{00}=\dots=\theta_{0,j-1}=0$, we
see from \eqref{Fr} that $F_0^{(r)}=\dots=F_{j-1}^{(r)}=0$ and
$F_{j}^{(r)}=\theta_{0j}$. It follows
$G_l^{(r)}(\be_r,\dots,\be_{r+rm+n-1})=0$ for $rm+n-j+1\le l\le rm+n$, $p\nmid l$,
and
\[
G_{rm+n-j}^{(r)}(\be_r,\dots,\be_{r+rm+n-1})=\be_r^{-m}\theta_{0j}\ne0.
\]
It follows
$h_\gp(L/K)=rm+n-j$ for any $\gp\in T_r$ by Lemma \ref{j}.

Let $p|rm+n-j$, then $h_\gp(L/K)=rm+n-j-1$ iff
$$
J_{r+rm+n-1}(\gp)=J_{r+rm+n-1}((-T+\la(\be_r)U^r+\dots+\la(\be_{r+rm+n-1})U^{r+rm+n-1})),
$$
where $\be_r,\dots,\be_{r+rm+n-1}$ are such that
$$
G_{rm+n-j-1}^{(r)}(\be_r,\dots,\be_{r+rm+n-1})
=\be_r^{-m}F_{j+1}^{(r)}(\be_r,\be_{r+1}/\be_r,\dots,\be_{r+j}/\be_r)\ne0.
$$
An easy computation shows that
$F_{j+1}^{(r)}(X_0,X_1,\dots,X_j)=\theta_{0,j+1}-m\theta_{0j}X_1$.
Thus, we have $h_\gp(L/K)=rm+n-j-1$ iff $\theta_{0,j+1}\be_r-m\theta_{0j}\be_{r+1}\ne0$.
Since $m\theta_{0j}\ne0$, such $\be_r,\dots,\be_{r+rm+n-1}$ exist.
We have proved that
$$
h_r(L/K)=\begin{cases} rm+n-j, & p\nmid rm+n-j,
         \\ rm+n-j-1, & p\mid rm+n-j.
         \end{cases}
$$
for $r>\max(1,j)$.

Case 2: $p\mid m$. We claim that there exists a positive integer $j$
such that $j\not\equiv n\mod p$, and $\theta_{0j}\ne0$. Indeed, assume that
$\theta_{0j}=0$ for all $j\not\equiv n\mod p$. Choose $N\ge n$ such that $p|N$.
Then the residue class of $T^mU^Na$ in $\widehat{A}/(T)=k[[U]]$ belongs to $k[[U^p]]$.
Consider the commutative diagram
$$
\CD
A_{(T)} & @>\al>> & {\widehat A}_{(T)}=k[[T,U]]_{(T)}
\\
@VVV & & @VVV
\\
\Fr(A/(T)) & @>\be>> & k((U))
\endCD
$$ where $\al$ is an embedding of discrete valuation rings and
$\be$ is the induced embedding of their residue fields. On the
other hand, $\be$ is exactly the completion map of the field of
fractions of the discrete valuation ring $A/(T)$. We can identify
$\Fr(A/(T))$ with the image of $\be$. Let $x\in k[[U]]$, $x^p\in
A/(T)$. Since $A$ is a G-ring, $A/(T)$ is also a G-ring, and any
algebraic subextension in $k((U))/\Fr(A/(T))$ is separable, whence
$x\in A/(T)$. It follows $$ k[[U^p]] \cap A/(T) = (A/(T))^p, $$
whence one can write $T^mU^Na=b^p+Tb'$, $b,b'\in A$. It follows
\[
v_1(a-T^{-m}U^{-N}b^p+T^{-m/p}U^{-N/p}b)>-m.
\]
 This yields a
contradiction with the minimality of $m$.

Now let $j$ be minimal with $\theta_{0j}\ne0$ and $j\not\equiv n\mod p$.
Changing $a$, we may assume without loss of generality
that $\theta_{0i}=0$ for all $i<j$.

Take any $r>\max(1,j)$. As in Case 1, we obtain
$h_\gp(L/K)=rm+n-j$ for all $\gp\in T_r$. (Note that $p\nmid rm+n-j$.)
 We have proved that  $h_r(L/K)=rm+n-j$ for all $r>\max(1,j)$.

In both cases we see that $\lim\limits_r h_r(L/K)/r=m$. \qed

\end{document}